# Nonparametric estimation of a distribution function under biased sampling and censoring

## Micha Mandel[1],*

*The Hebrew University of Jerusalem*

**Abstract:** This paper derives the nonparametric maximum likelihood estimator (NPMLE) of a distribution function from observations which are subject to both bias and censoring. The NPMLE is obtained by a simple EM algorithm which is an extension of the algorithm suggested by Vardi (*Biometrika*, 1989) for size biased data. Application of the algorithm to many models is discussed and a simulation study compares the estimator's performance to that of the product-limit estimator (PLE). An example demonstrates the utility of the NPMLE to data where the PLE is inappropriate.

## 1. Introduction

In this paper, the EM algorithm [4] of Vardi [18] is extended from size-biased to $W$-biased observations, where $W$ is a known, positive, increasing and right continuous function. Specifically, the algorithm finds the distribution function $G$ that maximizes

$$\prod_{i=1}^{m} \frac{dG(x_i)}{\mu^*} \times \prod_{j=1}^{n} \frac{\bar{G}(y_j)}{\mu^*}, \qquad (1.1)$$

where $\mu^* = \int_0^\infty W(x) dG(x)$, $\bar{G} = 1 - G$ and $x_1, \ldots, x_m$ and $y_1, \ldots, y_n$ are given data points.

The function (1.1) is proportional to likelihoods that arise in reliability and survival studies when data are subject to both bias and censoring. Several authors derive the nonparametric maximum likelihood estimator (NPMLE) of $G$ in problems where the likelihood is a special case of (1.1). The size-biased case, $W(x) = x$, appears in cross-sectional sampling if the population is in steady state. An early work is Cox [3] who estimates $G$ in the uncensored case (i.e., $n = 0$). Vardi [18] presents four problems that result in likelihood proportional to (1.1) with $W(x) = x$ and develops a simple EM algorithm to find the NPMLE of $G$. In an earlier paper [16], he develops an EM algorithm to estimate the underlying lifetime distribution of a renewal process observed in a time window (see also [15]). Wijers [22] suggests the EM algorithm for the same window sampling in a population model under

[1]Department of Statistics, The Hebrew University of Jerusalem, Mount Scopus, Jerusalem 91905, Israel, e-mail: `micha.mandel@huji.ac.il`
*This paper is part of my PhD dissertation that was completed under the supervision of Yosef Rinott and Yehuda Vardi. I am grateful for their excellent guidance and many suggestions that improved this paper considerably. A reviewer provided helpful comments and suggestions. Research supported by the Israel Science Foundation (grant No. 473/04).

the steady state assumptions. When the renewals or birth times are known for all sampled individuals, the likelihood is proportional to (1.1) with $W(x) = x + C$ where $C$ is the window width. Motivated by data on HIV infection, Kalbfleisch and Lawless [7] consider a population model in which entrances are according to an inhomogeneous Poisson process. Their likelihood is similar to the uncensored part of (1.1) with $W$ being the cumulative rate of the process. They briefly remark on the censored case at the end of their paper. The random left truncation model is commonly used for survival data [12, 21, 23]. In that model, two independent variables $A \sim W$ and $X \sim G$ are truncated to the region $A \leq X$. The variable $X$ may be censored by $A + C$, where $C$ is independent of $(A, X)$. When the truncation distribution is known, the likelihood of the model is proportional to (1.1). More models that give rise to likelihoods proportional to (1.1) are reviewed in Section 3.

In Section 2, we suggest a unified EM algorithm that provides the NPMLE for the general likelihood (1.1) and discuss its convergence properties. Two methods for developing the EM algorithm are considered; the first uses an extension of Vardi's multiplicative censoring model [18], and the second utilizes the random left truncation model. By presenting the two approaches, the similarity between the two seemingly unrelated models is highlighted. In Section 3 we derive the likelihoods of the above mentioned papers and other models and show that they are special cases of (1.1). In Section 4 we compare the performances of estimators with full and no knowledge on $W$ by simulation. The algorithm is then used to reanalyze the Channing House data [6]. We complete the paper with discussion in Section 5.

Throughout the paper, we will refer to (1.1) and to similar expressions as "likelihood", although they may be only proportional to the likelihood of the data. We will refer to the maximizer of (1.1) as the NPMLE of $G$. To distinguish, we will call the maximizer of $G$ when $W$ is not known the product-limit estimator (PLE). The latter is presented in Section 3.4 and is used in the simulation and application for comparison purposes.

## 2. An EM algorithm

There are several ways to develop the EM algorithm. We first adopt the approach used in [18] and solve a seemingly unrelated multiplicative censoring problem (similar to [18] problem A). This is shown to be almost equivalent to our original problem of maximizing (1.1). The second approach is motivated by cross-sectional sampling where subjects are selected to the sample at a random time point (similar to [18] problem B). It is assumed that ages at sampling are observed for all sampled individuals and residual lifetimes are subject to random censoring. Here the censoring times are independent, and are independent of the ages at sampling and the residual lifetimes. Both approaches eliminate the bias of the data and deal only with censoring. They demonstrate that bias is a secondary problem for estimation relative to censoring (when the bias is known and does not cause identification problems). After developing the algorithm, several convergence properties are discussed.

### 2.1. Vardi's approach

Since $W(x) > 0$ and known, (1.1) is proportional to

$$\text{(2.1)} \quad \prod_{i=1}^{m} dG^W(x_i) \prod_{j=1}^{n} \int_{y_j}^{\infty} \frac{1}{W(u)} dG^W(u)$$

where $G^W(x) = \int_0^x W(u)dG(u)/\mu^*$ is the $W$-weighted version of $G$. The problem can be divided into two parts; maximization of (2.1) for $G^W$, and transformation using $G(dx) \propto G^W(dx)/W(x)$. To maximize (2.1), consider first the following multiplicative censoring model which extends the problem studied by Vardi [18].

Let $X_1^0, \ldots, X_m^0, Z_1^0, \ldots, Z_n^0$ be positive random variables from the distribution $G^0$ supported on $(0, \infty)$ and let $U_1, \ldots, U_n \sim U(0,1)$ where all the random variables are independent. Let $Y_i^0 = W(Z_i^0)U_i$ $(i = 1, \ldots, n)$. Similar to Vardi who solves the problem for the special case $W(x) = x$, the $Y_i^0$'s describe the transformation from the 'complete data' $(X_1^0, \ldots, X_m^0, Z_1^0, \ldots, Z_n^0)$ to the observed 'incomplete data'. This transformation is used in the E-step of the EM algorithm. The statistical problem is of estimating $G^0$ using $(x_1^0, \ldots, x_m^0, y_1^0, \ldots, y_n^0)$, a realization of $(X_1^0, \ldots, X_m^0, Y_1^0, \ldots, Y_n^0)$.

First note that the density of $Y^0$ with respect to Lebesgue measure is

$$f_{Y^0}(t) = \int_{v \geq t} \frac{1}{v} dF_{W(Z^0)}(v) = \int_{\{v: W(v) \geq t\}} \frac{1}{W(v)} dG^0(v)$$

where $f_V$ and $F_V$ denote the density and distribution of a random variable $V$. The likelihood of the data $(\mathbf{x^0}, \mathbf{y^0}) = (x_1^0, \ldots, x_m^0, y_1^0, \ldots, y_n^0)$ is given by:

$$\mathcal{L}(G^0; \mathbf{x^0}, \mathbf{y^0}) = \prod_{i=1}^m dG^0(x_i^0) \prod_{j=1}^n \int_{\{v: W(v) \geq y_j^0\}} \frac{1}{W(v)} dG^0(v).$$

By defining $W^{-1}(y) = \min\{v : W(v) \geq y\}$ (which exists from the right continuity assumption) and recalling that $W$ is increasing, we can rewrite the likelihood as

$$(2.2) \qquad \mathcal{L}(G^0; \mathbf{x^0}, \mathbf{y^0}) = \prod_{i=1}^m dG^0(x_i^0) \prod_{j=1}^n \int_{v \geq W^{-1}(y_j^0)} \frac{1}{W(v)} dG^0(v)$$

and the similarity to our original likelihood (2.1) is apparent.

When $W$ is strictly increasing, no information is lost by the transformation $X^0 \mapsto W(X^0)$, and Vardi's EM algorithm [18] can be applied to $W(X_1^0), \ldots, W(X_m^0), Y_1^0, \ldots, Y_n^0$. This gives the NPMLE of $F_{W(X^0)}$ and the NPMLE of $G^0$ is obtained by a simple transformation. Specifically, let $t_1 < t_2 < \cdots < t_h$ be the distinct values of $x_1^0, \ldots, x_m^0, W^{-1}(y_1^0), \ldots, W^{-1}(y_n^0)$, and let $\xi_j, \zeta_j$ be the multiplicity of the $X^0$ and $W^{-1}(Y^0)$ samples at $t_j$:

$$\xi_j = \sum_{i=1}^m I\{x_i^0 = t_j\}, \qquad \zeta_j = \sum_{i=1}^n I\{W^{-1}(y_i^0) = t_j\}.$$

Denote by $p_j$ the mass $G^0$ assigns to $t_j$, $\mathbf{p} = (p_1, \ldots, p_h)$; then the problem can be rewritten as

$$(2.3) \qquad \begin{aligned} \text{maximize} \quad & \mathcal{L}(\mathbf{p}) = \prod_{j=1}^h p_j^{\xi_j} \left( \sum_{k \geq j} \frac{1}{W(t_k)} p_k \right)^{\zeta_j} \\ \text{subject to} \quad & \sum_j p_j = 1 \\ & p_j \geq 0 (j = 1, \ldots, h). \end{aligned}$$

And an EM step is

$$(2.4) \qquad p_j^{new} = (n+m)^{-1} \left( \xi_j + [W(t_j)]^{-1} p_j^{old} \sum_{k \leq j} \frac{\zeta_k}{\sum_{l \geq k} [W(t_l)]^{-1} p_l^{old}} \right),$$

where $p_j^{old}$ and $p_j^{new}$ are the current and updated estimates of $p_j$.

The derivation of (2.4) from (2.3) holds true whether or not $W$ is strictly increasing. Maximization of (2.2) reduces to the discrete problem (2.3) by showing that the support of the NPMLE is discrete and determining $t_1, \ldots, t_h$. This somewhat technical point is deferred to the Appendix for the case of a general non-decreasing $W$.

Finally, the connection between $G^0$ and $G^W$ is apparent from (2.1) and (2.2), where the only difference appears in the left limits of the integrals. These were used only to determine the points $t_1, \ldots, t_h$ and their multiplicity in the two samples. Thus, to use the algorithm for maximizing (2.1), one needs to define $t_1, \ldots, t_h$ as the distinct values of $x_1, \ldots, x_m, y_1, \ldots, y_n$ (the original data) and to change $\xi_j$ and $\zeta_j$ accordingly. The support of our original problem is a subset of the observations as discussed in the Appendix. After finding it, the problem reduces to (2.3) for which the algorithm (2.4) derives the NPMLE of $G^W$. The corresponding estimate of $G$ is achieved by the inversion formula $dG(x) \propto [W(x)]^{-1} dG^W(x)$, as mentioned above.

The Kaplan-Meier estimate of $G$ [8] is obtained by using $W \equiv 1$ in (2.4) (see Section 3.1 for more details). This estimator does not use the correct weights when it redistributes the mass of the censored observations, and in general it is inappropriate.

## 2.2. A direct approach

Cross-sectional samples of lifetimes usually contain the age at sampling $A$ and the residual lifetime $R$; the latter is subject to random censoring. Vardi's approach uses the sufficient statistic $A + R$ [see (2.7) below] to develop an EM algorithm. In this subsection, the statistic $(A, R)$ is used. The two approaches give different perspectives about the formation of the bias and censoring.

Assume that $W(0) = 0$, $W(t)\bar{G}(t) \to 0$ as $t \to \infty$, and $\int_0^\infty W(dt)\bar{G}(dt) = 0$. Mathematically this means that $\mu^* \equiv \int_0^\infty W(t)G(dt) = \int_0^\infty \bar{G}(t)W(dt)$. Practically it means that the probability of leaving the population at the very instant of sampling is zero; hence, one does not need to worry about inclusion or exclusion of such observations in the sample.

Let

$$(2.5) \qquad f_A(a) = \frac{\bar{G}(a)}{\mu^*} dW(a), \quad a > 0,$$

$$(2.6) \qquad f_{R|A}(r|a) = \frac{dG(a+r)}{\bar{G}(a)}, \quad r \geq 0$$

so the joint density at $(a, r)$ is

$$(2.7) \qquad f_{A,R}(a, r) = \frac{dG(a+r)}{\mu^*} dW(a), \quad a > 0, \ r \geq 0.$$

Now suppose we have one sample of $m$ pairs $(a_i, r_i)$ from $f_{A,R}$ and another independent sample of $n$ variables $y_j$ from $f_A$. This describes the so-called incomplete data and by denoting $x_i = a_i + r_i$ $(i = 1, \ldots, m)$ we arrive at the likelihood (1.1). Since $W$ is known, the product $\prod_i dW(a_i) \prod_j dW(y_j)$ is irrelevant for maximization.

The complete data are of course $x_1, \ldots, x_m$ and $y_1 + \tilde{r}_1, \ldots, y_n + \tilde{r}_n$, where $\tilde{r}_j$ is the unobserved residual lifetime of subject $j$ of the second sample. Using the sum $x_i = a_i + r_i$ instead of its components is justified by noticing that the sum is the sufficient statistic for the complete data problem. By variables changing and integrating $a$ out in (2.7), it can be easily verified that the likelihood of $x_i$ (or the density of $y_j + \tilde{r}_j$) is

$$(2.8) \qquad \frac{W(x_i)dG(x_i)}{\mu^*}.$$

For the support points $t_1, \ldots, t_h$ described in the appendix, the E-step uses (2.6):

$$(2.9) \qquad E_{\mathbf{\Pi}^{old}}\left(I\{A_i + R_i = t_j\}|A_i = y_i\right) = \frac{\pi_j^{old}}{\sum_k \pi_k^{old} I\{t_k \geq y_i\}} I\{t_j \geq y_i\},$$

where $\mathbf{\Pi}^{old} = (\pi_1^{old}, \ldots, \pi_h^{old})$ is the current estimate of the unbiased distribution $G$, i.e., the estimate at $t_j$ of the weighted distribution $G^W$ is $p_j^{old} \propto W(t_j)\pi_j^{old}$.

The complete likelihood is a product of terms such as (2.8) which is the likelihood of a weighted sample. An M-step estimates the weighted distribution $G^W$ by the empirical distribution function. Combining the E-step and the M-step, an iteration is given by:

$$(2.10) \qquad p_j^{new} = (n+m)^{-1}\left(\xi_j + \sum_{k \leq j} \frac{\zeta_k \pi_j^{old}}{\sum_{l \geq k} \pi_l^{old}}\right),$$

where $\xi_j$ and $\zeta_j$ are the multiplicities of uncensored and censored observations, respectively. Put $\pi_j^{old} = [W(t_j)]^{-1} p_j^{old} / \sum_k \{[W(t_k)]^{-1} p_k^{old}\}$ in the equation above to get (2.4).

### 2.3. Convergence of the algorithm

Several properties of the problem and the algorithms are sketched below. The appendix shows that for some functions $W$, such as step functions, the NPMLE is not unique and different choices of the support points yield the same (maximum) value of the likelihood. The properties below hold for a given choice of support points.

**Property 2.1.** *Given the points $t_1, \ldots, t_h$, the maximizer of (2.3) is unique.*

*Proof.* We show that the problem can be replaced with a maximization of a strictly concave function over a convex region. Following [18], write $q_j = p_j/W(t_j)$, $Q_j = \sum_{k \geq j} q_k$. Since $W(t_j)$ is constant in the likelihood, we can replace (2.3) with:

$$(2.11) \qquad \begin{aligned} \text{maximize} \quad & \log\Big(\prod_{j=1}^{h} q_j^{\xi_j} Q_j^{\zeta_j}\Big) \\ \text{subject to} \quad & \sum_j q_j W(t_j) \leq 1 \\ & q_j \geq 0 \, (j = 1, \ldots, h). \end{aligned}$$

Now $\log \left( \prod_{j=1}^{h} q_j^{\xi_j} Q_j^{\zeta_j} \right) = \sum_{j=1}^{h} \xi_j \log(q_j) + \sum_{j=1}^{h} \zeta_j \log(Q_j)$, and the assertion follows from log of a partial sum being a strictly concave function from $\mathbb{R}^{+h}$ to $\mathbb{R}$, and since sum of concave functions is concave (recall that for all $j$, $\xi_j, \zeta_j \geq 0$ and $\xi_j + \zeta_j \geq 1$). □

**Property 2.2.** *Let* $\mathbf{p}_n$ *and* $\mathcal{L}_n = \mathcal{L}(\mathbf{p}_n)$ *be the value of* $\mathbf{p}$ *and the value of the likelihood assigned by the EM algorithm (2.4) after its n'th iteration; then* $\mathcal{L}_n$ *converges to a point* $\mathcal{L}^*$.

*Proof.* [1] and later [4] show that the likelihood increases in each iteration of the EM algorithm. The assertion follows from (2.3) that shows that the likelihood is bounded from above. □

**Property 2.3.** *If the maximizer of (2.3) assigns positive mass to all points* $t_1, \ldots, t_h$, *then the algorithm (2.4) converges to that maximizer.*

*Proof.* This follows from Theorem 4 of [24] by using the uniqueness proved in Property 2.1 and the fact that $\mathcal{L}$ is a polynomial on the simplex (which establishes the regularity conditions). □

## 3. Examples

Many examples are easily described using the Lexis diagram (see Figure 1) that depicts changes in a population $\mathcal{S}$ over time. The horizontal and vertical axes of the diagram represent the calendar time and the lifetime or duration of subjects in $\mathcal{S}$, respectively. Subjects of $\mathcal{S}$ are represented by 45° lines that start at their time of entering $\mathcal{S}$ and end at leaving. Lines that cross the vertical line $x = t$ correspond to the population of $\mathcal{S}$ at $t$. In particular, a cross-sectional sample at time 0 contain those subjects (or lines) that intersect with the line $x = 0$. For a review of the Lexis diagram and its utility for studying different sampling designs and population quantities see [2] and [11].

### *3.1. Random censorship*

A common sampling plan is of collecting data on subjects entering $\mathcal{S}$ during the time window $[0, C]$. This widely used design is known as the random censorship model and it is depicted in the top left panel of Figure 1. According to the model, a random sample $T_i$ $(i = 1, 2, \ldots, m + n)$ is selected from a distribution $G$, but instead of the $T_i$'s, observations are independent realizations of $\min(T_i, C_i)$, where $C_i \sim F_C$ is independent of $T_i$. The data contain also the information whether the observations were censored or not.

Denote the uncensored observation by $x_1, \ldots, x_m$ and the censored ones by $y_1, \ldots, y_n$, then the likelihood of the data is

$$(3.1) \qquad \prod_{i=1}^{m} dG(x_i) \prod_{j=1}^{n} \bar{G}(y_j) \times \prod_{i=1}^{m} \bar{F}_C(x_i) \prod_{j=1}^{n} dF_C(y_j).$$

The NPMLE is the celebrated Kaplan-Meier estimator [8] and is based only on the left part of (3.1). This is a special case of (1.1) with $W = 1$, and can be solved (in somewhat a redundant way) by the EM algorithm (2.4). At convergence, the algorithm is exactly the self-consistent estimate of Efron [5] and (2.4) illustrates the redistributed to the right principle. As Efron shows, it reduces to the non-iterative Kaplan-Meier estimator.

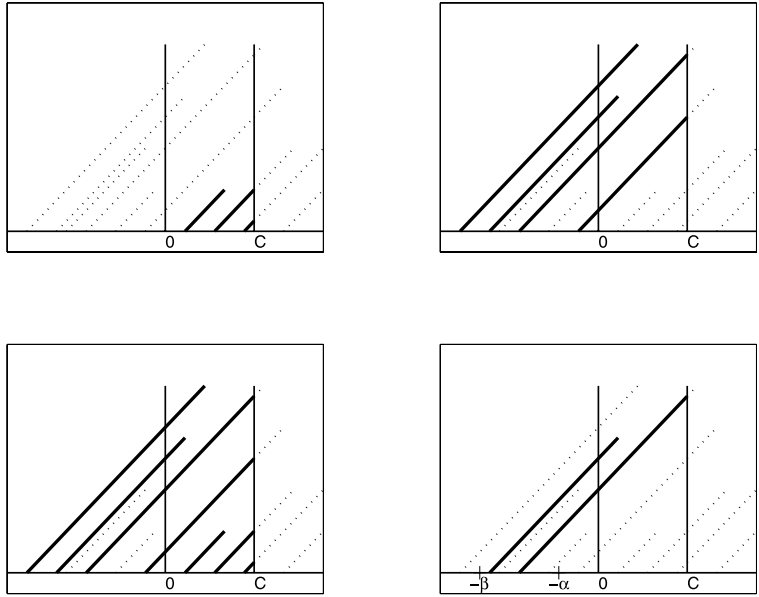

FIG 1. *Sampling designs in the Lexis diagram. Each $45°$ line represents a subject in the population. The horizontal axis is the calendar time and the vertical axis is duration or lifetime. Observed data are depicted as solid lines.*
*Top left - random censoring model. $\mathcal{I} = (0, C)$*
*Top right - follow-up on cross-sectional sampling. $\mathcal{I} = (-\infty, 0)$.*
*Bottom left - window sampling. $\mathcal{I} = (-\infty, C)$.*
*Bottom right - cross-sectional with truncation. $\mathcal{I} = (-\beta, -\alpha)$.*

According to Lemma A.3, an NPMLE assigns mass only to the uncensored observations, because $W$ is constant everywhere (if the last observation is censored, then it assigns mass also to $\max\{y_j\}$). This is a well-known feature of the Kaplan-Meier estimator.

### 3.2. Population models with poisson entrance process

The likelihood (1.1) is obtained in many designs where entrances to the studied population $\mathcal{S}$ are according to a Poisson process. Specifically, the model assumes a Poisson entrance process $N(t)$ on $(-\infty, C)$ with rate $\rho(\cdot)$. The lifetimes, $X_1, X_2, \ldots$, are determined by the law $G$, and $N(\cdot), X_1, X_2, \ldots$ are independent. The sample consists of all subjects who entered $\mathcal{S}$ during $\mathcal{I}$ and are in $\mathcal{S}$ sometime during $(0, C)$, where $\mathcal{I}$ is a subset of $(-\infty, C)$, usually an interval. Figure 1 presents several common examples and the corresponding samples. It is assumed that $\rho(x) = \lambda \times \lambda_0(x)$ where $\lambda_0$ is known and that $\int_{\mathcal{I}} \lambda_0(u)\bar{G}(-u)du = \mu^* < \infty$. For each sampled subject, we observe the entrance time $a$ and the possibly censored residual lifetime $r$ with the failure indicator $\delta$. Denote by $(a_1, x_1), \ldots, (a_m, x_m)$ and $(a_{m+1}, y_{m+1}), \ldots, (a_{m+n}, y_{m+n})$ the data on the $m$ uncensored and $n$ censored observations, then the likelihood is

given by

$$(3.2) \quad \mathcal{L}(G) = \prod_{i=1}^{m} \frac{dG(x_i)\lambda_0(a_i)}{\mu^*} \prod_{j=m+1}^{m+n} \frac{\bar{G}(y_j)\lambda_0(a_j)}{\mu^*} \times e^{-\lambda\mu^*} \frac{(\lambda\mu^*)^{n+m}}{(n+m)!}.$$

The NPMLE of $\lambda\mu^*$ is $n + m$ and since $\lambda_0$ is assumed known, the problem reduces to maximizing for $G$ the function

$$(3.3) \quad \mathcal{L}(G|n+m) = \prod_{i=1}^{m} \frac{dG(x_i)}{\mu^*} \prod_{j=m+1}^{m+n} \frac{\bar{G}(y_j)}{\mu^*}$$

which has the form of (1.1). In this case $W(x) = \int_{(-x,C)\cap \mathcal{I}} \lambda_0(u)du$ which is continuous and increasing.

**Examples**:

*A homogeneous poisson process - cross-sectional sampling* - $\mathcal{I} = (-\infty, 0)$. Early studies of this design are [14] and [9]. It is depicted in Figure 1 in the top right panel. The sample consists of all subjects who are in $\mathcal{S}$ at time 0 and the data contain their entrance time and follow-up until time $C$. Here $W(x) = x$ and $\mu^* = \mathbb{E}X$ is the mean lifetime. The algorithm (2.4) for this special case is derived by [18] who studies a slightly different model.

*A homogeneous poisson process - window sampling* - $\mathcal{I} = (-\infty, C)$. Here the sample consists of the cross-sectional population and those entering during the follow-up period (see the bottom left panel of Figure 1). Thus, there are two samples: i) a size biased sample comprising of subjects who entered $S$ before 0, and ii) an unbiased sample comprising of those who entered during the time window $[0, C]$. Both samples are censored at $C$. This model is a mixture of the random censorship and the cross-sectional models discussed above, with random number of observations from each model. Here $W(x) = x + C$ and $\mu^* = \mu + C$. The model is studied by [9]; [22] provides an EM algorithm for it.

*An inhomogeneous poisson process - a window sampling* - $\mathcal{I} = (-\infty, C)$. Kalbfleisch and Lawless [7] study entrances according to an inhomogeneous process in a somewhat different model, mainly focusing on the uncensored case. They consider many models in this framework and estimate both the rate function and the distribution of lifetimes. Under this model, $W(x) = \int_{-x}^{C} \lambda_0(t)dt$ is proportional to the cumulative rate.

*A truncated poisson process* - $\mathcal{I} = (-\beta, -\alpha)$. Wang [20] derives the NPMLE of $G$ when data are collected only for subjects whose ages at sampling time are in $[\alpha, \beta]$ for some known $0 \leq \alpha < \beta$. This design is depicted in the bottom right panel of Figure 1 and it is natural when data started to be recorded only $\beta$ years ago (in that case $\alpha = 0$) or when there is a specific interest on subjects who entered during the period $(-\beta, -\alpha)$ (e.g., when $\mathcal{S}$ is defined by some epidemic status and a specific treatment was used during $(-\beta, -\alpha)$). Wang shows that when entrances to $\mathcal{S}$ are according to a homogeneous Poisson process the likelihood is given by

$$\prod_i \frac{dG(x_i)}{\int_\alpha^\beta \bar{G}(u)du} \prod_j \frac{\bar{G}(y_j)}{\int_\alpha^\beta \bar{G}(u)du}.$$

A simple calculation shows that the likelihood is of the form (1.1) with $\mu^* = \mathbb{E}[\min(X, \beta) - \alpha]^+$ so the algorithm (2.4) can be applied with $W(x) = [\min(x, \beta) - \alpha]^+$. We note that in this case $W(x) = 0$ for $x \in [0, \alpha]$ so implementation of the algorithm looks problematic. However, under this setting, $G$ is not identifiable on $[0, \alpha]$ [20] and one can only hope to estimate $G$ given $X > \alpha$ where $W > 0$. This model is used in Section 4 to analyze the Channing House data [6].

### 3.3. A discrete entrance process

Mandel and Rinott (unpublished) study a discrete version of the Poisson entrance model in which entrances to $\mathcal{S}$ occur at fixed time points $\sigma_K < \cdots < \sigma_2 < \sigma_1 \leq 0$. At $\sigma_k$, $N_k$ new subjects joining $\mathcal{S}$ with lifetimes $X_{k1}, \ldots, X_{kN_k}$. The model assumes that $N_k$ has a Poisson distribution with parameter $\lambda \times \lambda_0(k)$, where $\lambda_0$ is known, $X_{ki} \sim G$ ($k = 1, \ldots, K$; $i = 1, \ldots, N_k$), and $\{N_1, \ldots, N_K, X_{11}, X_{12}, \ldots, X_{KN_K}\}$ are independent. Under this model, the NPMLE is obtained by maximizing (1.1) with
$$W(x) = \sum_{\{k: -\sigma_k \leq x\}} \lambda_0(k).$$
Here $W$ is a step function and Lemma A.3 should be used to determine the support of $G$. When $\sigma_k = -k$, $\lambda_0$ is constant and lifetimes are integer valued, the model is a discrete size-biased model.

### 3.4. Truncation models

The left truncation model [12, 21, 23] assumes that two independent variables $A \sim W$ and $T \sim G$ can be observed only on the region $A \leq T$. In addition, there is a random censoring variable $C$ which is independent of $T$ and satisfies $P(C > A) = 1$. Data comprise of $(a_i, \min(t_i, c_i), \delta_i)$ ($i = 1, \ldots, m+n$) which are $m+n$ realizations of $(A, \min(T, C), \Delta)$ restricted to the region $A \leq T$, where $\Delta = 1$ if $T \leq C$ and $\Delta = 0$ otherwise. Let $\mu^* = P(A \leq T) = \mathbb{E}W(T)$. Changing the notations as in (3.2), the likelihood of the data is

$$(3.4) \qquad \prod_{i=1}^{m} \frac{dW(a_i)dG(x_i)}{\mu^*} \prod_{j=m+1}^{m+n} \frac{dW(a_j)\bar{G}(y_j)}{\mu^*}$$

(note the similarity to (3.2)). If $W$ is known, then $dW(u_i)$ can be omitted from the likelihood and the problem of maximizing (3.4) is equivalent to maximizing (1.1).

Equation (3.4) can be reexpressed as

$$(3.5) \quad \prod_{i=1}^{m} \frac{dG(x_i)}{\bar{G}(a_i-)} \prod_{j=m+1}^{m+n} \frac{\bar{G}(y_j)}{\bar{G}(a_j-)} \times \prod_{i=1}^{m} \frac{dW(a_i)\bar{G}(a_i-)}{\mu^*} \prod_{j=m+1}^{m+n} \frac{dW(a_j)\bar{G}(a_j-)}{\mu^*},$$

where the first term is the likelihood of $T|A = a, A \leq T$ and the second term is the likelihood of $A|A \leq T$. If $W$ is completely unknown, maximizing (3.4) for $G$ is equivalent to maximizing the left term in (3.5) and the maximizer is the product-limit estimator (PLE) defined by

$$(3.6) \qquad \frac{\tilde{G}(dt)}{1 - \tilde{G}(t-)} = \frac{\sum_{i=1}^{m} I\{x_i = t\}}{\sum_{i=1}^{m} I\{a_i \leq t \leq x_i\} + \sum_{j=m+1}^{m+n} I\{a_j \leq t \leq y_j\}},$$

see [20]. In the next section, the PLE and the NPMLE are compared on real and simulated data.

## 4. Illustration

*Channing House data.* The data set comprises of 97 males who were residents of the Channing House retirement community in Palo Alto, California [6]. It contains the age at entry to the community and the age at death or censoring accompanied by the event indicator. In terms of Section 3.4, the age at entry is $A$ and the age at death is $T$. The parameter of interest, $G$, is the distribution of age at death (in months) of male residents of the community.

Wang [20] estimates the distribution of age (in months) at entry to be uniform on $(782, 1073)$. Assuming $W$ is the uniform distribution function on $(782, 1073)$, $G(\cdot)/\bar{G}(782)$ was estimated using algorithm (2.4) and $W(x) = [\min(x, 1073) - 782]^+$. Note that this is similar to the bias of the truncated Poisson process model described at the end of Section 3.2. The NPMLE and the PLE are depicted in the right panel of Figure 2. The survival under the uniform assumption is estimated to be somewhat higher.

For the analysis in the right panel of Figure 2, only 93 out of the 97 residents were used. Using all 97 individuals was problematic since the PLE approaches zero after the first two failures. This is shown in the left panel of Figure 2. Also shown is the survival estimated by (2.4) based on all 97 individuals and after estimating $W$ to be $U(751, 1073)$ (by the minimum and maximum age at entry of all 97 residents). The NPMLE is seen to be less sensitive to outliers. It shows that (2.4) can provide reasonable estimates when the PLE fails, a common phenomenon in small truncated data sets.

*Simulation.* A simulation study was conducted to compare the performances of the NPMLE and the PLE under the left truncation model described in Section 3.4. We used the EXP(1) model for both $G$ and $W$ and generated 400 data sets of 50 observations each. The first row of Figure 3 compares the performance of the estimators in terms of log MSE at the deciles of $G$, calculated by the average over the repeated replicates. Since the PLE is not well defined when the risk group is empty before the last observation [21], such data sets were not used to calculate the MSE of the PLE. They were used to calculate the MSE of the NPMLE. The columns of Figure 3 show the effect of censoring. Lifetime were censored at $A + C$ for a fixed $C$ such that the probability of censoring is 10, 25 and 50 percent from left to right. The second row of Figure 3 shows the results of the same analysis applied to data sets of 200 observations. These were generated by combining the simulated data sets of 50 observations (total of 100 data sets). Figure 4 shows the sensitivity of the estimators to the assumption on $W$. The analysis was done assuming the exponential distribution as before while the data were generated using $W=\text{Gamma}(2,1)$.

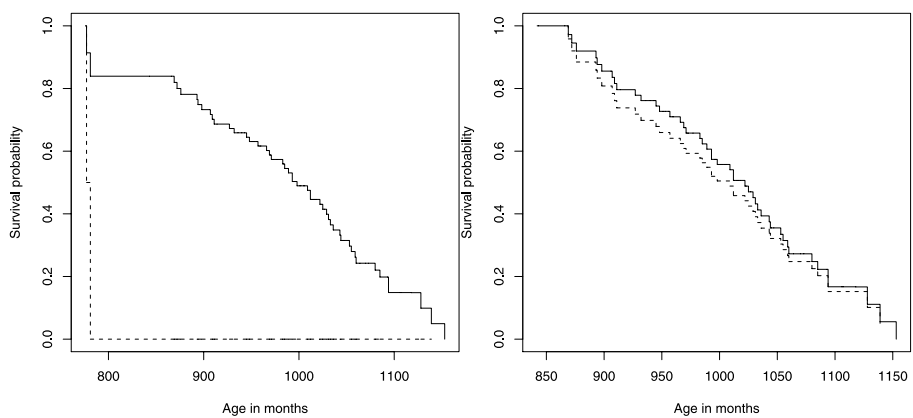

FIG 2. *Comparison of the PLE (dashed line) and the NPMLE (solid line) of survival in the Channing House community. Left - all 97 individuals. Right - excluding the first two failures.*

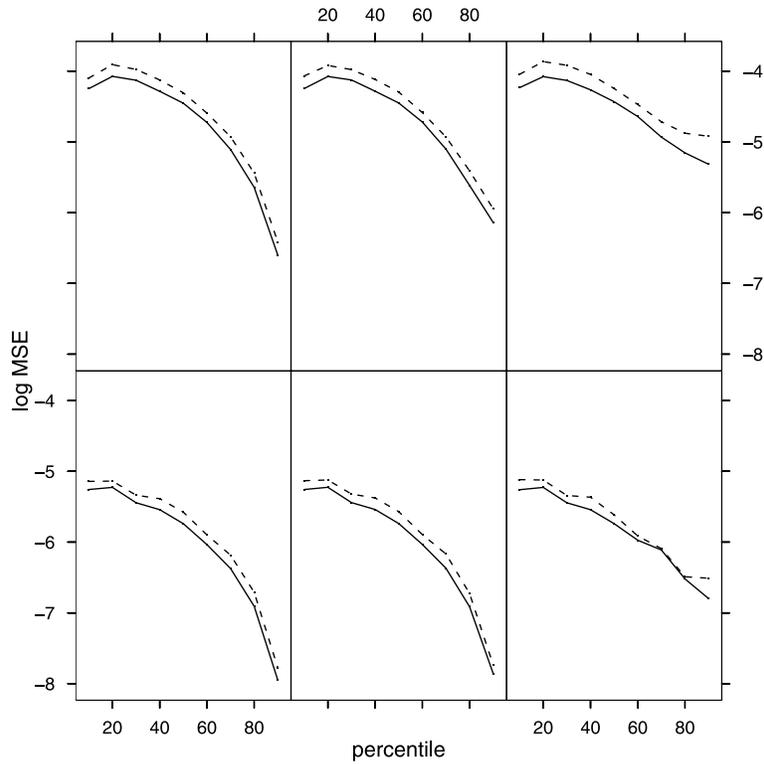

Fig 3. *Comparison of log MSE of the NPMLE (solid line) and the PLE (dashed line) calculated at deciles. Data were generated using $W=EXP(1)$. Sample size of 50 (top) and 200 (bottom), censoring probability of 0.1, 0.25 and 0.5 from left to right.*

The results show that using knowledge on the bias improves estimation by 10%-25% in terms of MSE. Increasing the probability of censoring results in a higher MSE of both estimators especially in the right tail. The relative performance does not change much by censoring, but some indication of better relative performance of the NPMLE in the right tail is seen in the simulation with 50% probability of censoring. More interesting are the results of the sensitivity study that show that the NPMLE is quite sensitive to the assumption on $W$. In general, the MSE of the PLE is smaller than that of the NPMLE. Moreover, the performance of the PLE improves when sample size increases from 50 to 200 while the performance of the NPMLE does not change. However, the performance of the NPMLE is better than that of the PLE in the left tail even when the model is incorrect. This phenomenon is seen in simulation with other distributions (not shown) and is probably attributed to the small risk group in the left tail that results in unstable estimation.

In the simulated data sets, the PLE was not well defined in 2% to 20% of the samples depending on the setting.

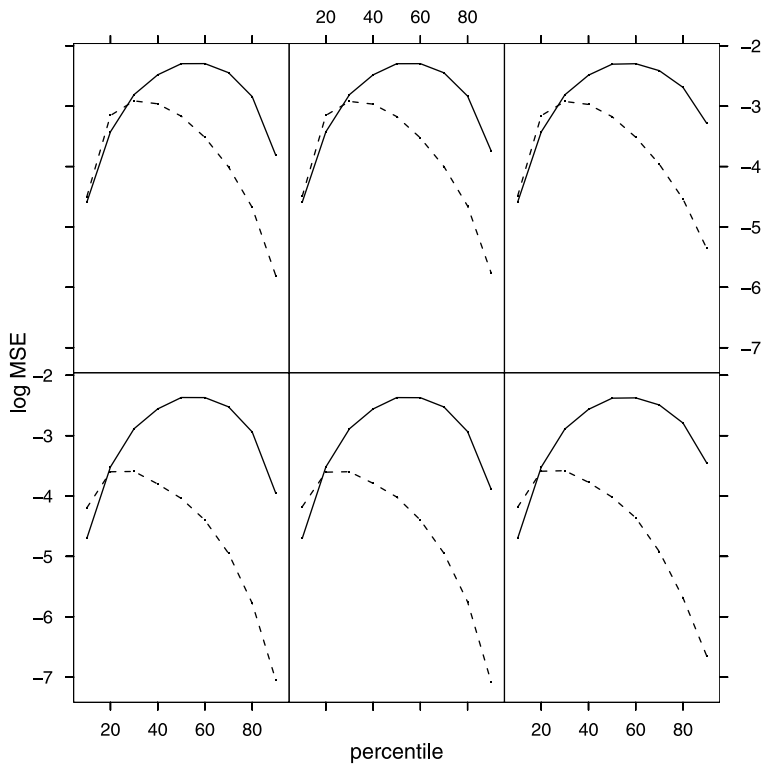

FIG 4. *Sensitivity analysis. Log MSE of the NPMLE (solid line) and the PLE (dashed line) at deciles. Data were generated using W=Gamma(2,1) and the NPMLE was calculated assuming W=EXP(1). Sample size of 50 (top) and 200 (bottom), censoring probability of 0.1, 0.25 and 0.5 from left to right.*

## 5. Discussion

The principal aim of this article is to provide a general framework and a unified algorithm for problems involving bias and censoring. A secondary aim is to contrast Vardi's multiplicative censoring model with truncated data and to compare the NPMLE to the PLE. An important question is which of the two estimators to use. The PLE cannot be calculated when all observations are censored or when data do not contain the truncation times $a_1, \ldots, a_{m+n}$ (see Section 3.4), and may not be well defined in data sets that do contain the truncation times as illustrated by the Channing House example. In all of these situations, the NPMLE exists and can be used. The simulation study shows that the NPMLE is more efficient when the model is correctly specified. However, it also indicates that it is quite sensitive to the assumed form of $W$. The use of the NPMLE, therefore, should be limited to situations where there is a theoretical justification for the assumed model of $W$. Furthermore, when data contain truncation times, the assumed form of $W$ can and should be tested. Wang [20] suggests a graphical goodness-of-fit test by plotting an estimate of $W$ versus the assumed model, and [10] study generalized Pearson statistics that can provide formal goodness-of-fit tests for the current model. These and other goodness-of-fit tests are studied in [13].

Algorithm (2.4) can be nested in more complex algorithms to provide nonpara-metric estimates for other interesting problems. Several examples are given in the unpublished PhD dissertation of the author. For example, Wang [19] studies the semi-parametric left truncation model $\{W \in \mathcal{W}_\theta, G \text{ unrestricted}\}$, where $\mathcal{W}_\theta$ is a family of distributions indexed by $\theta$. Her method, however, is applicable only for uncensored data and hence is of limited use. Estimates under Wang's model for censored data can be obtained by an iterative algorithm that uses (2.4) in one of its steps. Preliminary simulation results reveal better performance of this estimator over the PLE.

The algorithm presented in this article can be easily extended to likelihood of the form

$$(5.1) \qquad \prod_{s=1}^{S} \Big\{ \prod_{i=1}^{m_s} \frac{dG(x_{si})}{\mu_s^*} \times \prod_{j=1}^{n_s} \frac{\bar{G}(y_{sj})}{\mu_s^*} \Big\},$$

where $\mu_s^* = \int_0^\infty W_s(x) dG(x)$ for known increasing and right continuous functions $W_s$ ($s = 1, \ldots, S$). This likelihood generalized the model of [17] to the multiplicative censoring case. The complete likelihood in this problem involves products of likelihoods from different weight functions. The E-step is equivalent to (2.9), and the M-step uses Vardi's algorithm for selection bias models [17].

## Appendix A: The support of the NPMLE

This appendix discusses the determination of the support of the NPMLE for a non increasing right continuous function $W$. Although the NPMLE is not always unique (when $W$ has steps) it is shown that there exists an NPMLE that assigns mass only to the observed points (Lemma A.1). Furthermore, Lemmas A.2 and A.3 characterizes observed points that can be excluded from the support.

Let $W^{-1}(y) = \min\{v : W(v) \geq y\}$, then for the likelihood (2.2) we have

**Lemma A.1.** *There exists an NPMLE of $G^0$ that assigns mass only to the critical points $x_1, \ldots, x_m, W^{-1}(y_1), \ldots, W^{-1}(y_n)$.*

**Proof.** If mass is assigned to points other than the critical ones (i.e., points other than $x_j$ or $W^{-1}(y_i)$), then shifting the mass to the closest critical point to the left will not decrease the likelihood. The assertion follows after noticing that mass assigned to the left of the minimal critical point contributes nothing to the likelihood since the integrals have left limits.

Next, suppose that $m = n = 1$, $W$ is constant on an interval $[a, b)$ and $a < x_1 < b$ and $W^{-1}(y_1) = a$. Drawing a step function $W$ and inspecting the likelihood show that the NPMLE assigns mass only to $x_1$. In general, if there exist $i$ and $j$ such that $W(x_i) = W(W^{-1}(y_j))$, then the likelihood (2.2) increases if mass first assigned to $W^{-1}(y_j)$ is shifted to $x_i$. This excludes several of the critical points from being in the support of the NPMLE of $G^0$ and it is summarized by

**Lemma A.2.** *Let $W^{-1}(y|\mathbf{x}) = \min_{1 \leq i \leq m}\{x_i : W(x_i) = W(W^{-1}(y))\}$ if such $x_i$ exists, and $W^{-1}(y|\mathbf{x}) = W^{-1}(y)$ otherwise. Then there exists an NPMLE of (2.2) that assigns mass only to the points $x_1, \ldots, x_m, W^{-1}(y_1|\mathbf{x}), \ldots, W^{-1}(y_n|\mathbf{x})$.*

The support of the NPMLE of $G^W$ from the likelihood (2.1) is determined by arguments similar to those leading to Lemma A.1 and Lemma A.2. This suggests that the support is a subset of $\{x_1, \ldots, x_m, y_1, \ldots, y_n\}$ (from the integrals appearing in the likelihood (2.1) it is seen that the support points are the $y$'s and not the $W^{-1}(y)$'s). Moreover, if observations $y_j < x_i$ exist such that $W(x_i) = W(y_j)$, then the likelihood increases if mass initially assigned to $y_j$ is shifted to $x_i$. Likewise, if there are $y_j < y_{j'}$ such that $W(y_j) = W(y_{j'})$ then the likelihood increases if mass initially assigned to $y_j$ is shifted to $y_{j'}$. The following lemma summarizes this discussion:

**Lemma A.3.** *There always exists an NPMLE for (1.1) which assigns mass only to observed points. The complete observations $x_i$ $(i = 1, \ldots, m)$ are always points of support. A censored observation $y_j$ is not a point of support if: (i) there exists $x_i > y_j$ such that $W(x_i) = W(y_j)$ or (ii) there exists $y_{j'} > y_j$ such that $W(y_{j'}) = W(y_j)$.*

**Remark A.1.** Lemma A.3 tells us which $y_j$'s are not points of support and not which $y_j$'s are. The EM algorithm may assign mass zero to some of the $y_j$'s. One can always use the inefficient approach of considering all observations as support points and let the algorithm to assign zero mass to the redundant ones.